\documentclass[10pt,reqno]{amsart}
\usepackage[utf8]{inputenc} % allow utf-8 input
\usepackage{url, booktabs, amsfonts, nicefrac, microtype}      % microtypography

\usepackage[english]{babel} % Swedish translations
\usepackage{calc, amsmath, amssymb, amsthm}
\usepackage[english]{babel}
\usepackage[utf8]{inputenc}

\newcommand{\CC}{\mathbb{C}}
\newcommand{\QQ}{\mathbb{Q}}
\newcommand{\NN}{\mathbb{N}}

\newcommand{\mcA}{\mathcal{A}}
\newcommand{\mcC}{\mathcal{C}}
\newcommand{\mcS}{\mathcal{S}}

\newcommand{\mcO}{\mathcal{O}}

\newcommand{\mcE}{\mathcal{E}}

\newcommand{\mfP}{\mathfrak{P}}

\newcommand{\Li}{\mbox{Li}}

\newtheorem{thm}{Theorem}

\newtheorem{lemma}[thm]{Lemma}

\newcommand{\Rea}{\mathfrak{Re}}
\newcommand{\Ima}{\mathfrak{Im}}

\usepackage{import}

\usepackage{xcolor}

\begin{document}

\title{On an explicit zero-free region for the Dedekind zeta-function}
\author{Ethan~S.~Lee}
\address{UNSW Canberra, Northcott Drive, Campbell, Canberra, ACT 2612 Australia}

\maketitle

\begin{abstract}
%% Text of abstract
We establish new explicit zero-free regions for the Dedekind zeta-function. Two key elements of our proof are a non-negative, even, trigonometric polynomial and explicit upper bounds for the explicit formula of the so-called differenced logarithmic derivative of the Dedekind zeta-function. The improvements we establish over the last result of this kind come from two sources. First, our computations use a polynomial which has been optimised by simulated annealing for a similar problem. Second, we establish sharper upper bounds for the aforementioned explicit formula.
\end{abstract}

\section{Introduction}

Let $K$ be an algebraic number field and $L$ be a normal extension of $K$ with Galois group $G=\mbox{Gal}(L/K)$. Suppose $d_L$, $d_K$ denote the absolute values of the respective discriminant, $n_L = [L:\QQ]$ and $n_K = [K:\QQ]$. The Dedekind zeta-function of $L$ is denoted and defined for $\Rea(s)>1$ by
$$\zeta_L(s)= \sum_{\mfP}\frac{1}{N(\mfP)^s},$$
where $\mfP$ ranges over the non-zero ideals of $\mcO_L$. If $n_L = a + b$, then one can also consider the completed zeta-function
\begin{align*}
    \xi_L(s)&=s(s-1){d_L}^{\frac{s}{2}}\gamma_L(s)\zeta_L(s)\mbox{ such that}\\
    \gamma_L(s)&= \pi^{-\frac{as}{2}} \Gamma\left(\frac{s}{2}\right)^a \pi^{-\frac{b(s+1)}{2}} \Gamma\left(\frac{s+1}{2}\right)^b.
\end{align*}
Here, $\xi_L$ is an entire function satisfying the functional equation $\xi_L(s)=\xi_L(1-s)$. It can be seen that $\zeta_L$ is meromorphic on the complex plane with exactly one simple pole at $s=1$. Let $\mathcal{P}$ denote a prime ideal of $K$ and $P$ denote a prime ideal of $L$. If $\mathcal{P}$ is unramified in $L$, then the Artin symbol, $$\left[\frac{L/K}{\mathcal{P}}\right],$$ denotes the conjugacy class of Frobenius automorphisms corresponding to prime ideals $P|\mathcal{P}$. For each conjugacy class $C\subset G$, the prime ideal counting function is
$$\pi_C(x, L/K) = \#\left\{\mathcal{P} : \mathcal{P}\mbox{ unramified in }L, \left[\frac{L/K}{\mathcal{P}}\right] = C, N_K(\mathcal{P})\leq x\right\}.$$
In 1926, Chebotar\"{e}v \cite{chebotarev1926} proved the Chebotar\"{e}v density theorem, which states that
$$\pi_C(x,L/K)\sim \frac{\# C}{\# G}\Li(x) = \frac{\# C}{\# G}\int_2^x\frac{dt}{\log t}\mbox{ as }x\rightarrow\infty.$$
For example, if $L=K=\QQ$, then the Chebotar\"{e}v density theorem restates the prime number theorem. Moreover, if $\omega_\ell = e^{\frac{2\pi i}{\ell}}$ is the $\ell$th root of unity, $K=\QQ$ and $L=\QQ(\omega_\ell)$, then the Chebotar\"{e}v density theorem identifies with the Dirichlet theorem for primes in arithmetic progressions.

In 1977, Lagarias--Odlyzko \cite{lagarias1977effective} provided explicit estimates for the error term of the Chebotar\"{e}v density theorem. There are two results contained therein; one version assumes the generalised Riemann hypothesis (GRH) for $\zeta_L$ and the other does not. Their error term is effectively computable, dependent only on $x$, $n_L$, $d_L$ and $\frac{\# C}{\# G}$.

Under the GRH for $\zeta_L$, one can obtain the best possible effective results. Without assuming the GRH for $\zeta_L$, the better the zero-free region for $\zeta_L$ one has, the better the effective result one can achieve.
Therefore, the objective of this paper is to improve the best known, explicit zero-free region for $\zeta_L$, given by Kadiri \cite{kadiri2012explicit} in 2012.
We recall two famous forms of zero-free regions for the Riemann zeta-function.

\textit{Classical zero-free region.}
In 1899, de la Vall{\'e}e Poussin \cite{ValeePoussin} famously proved that there exists a positive constant $R$ such that $\zeta$ is non-zero in the region $s=\sigma + it$ such that $t\geq T$ and
\begin{equation}\label{eqn:dlvp}
    \sigma \geq 1-\frac{1}{R\log t}.
\end{equation}
The best known zero-free region for $\zeta$ of this kind is attributed to Mossinghoff--Trudgian \cite{MossinghoffTrudgian2015}, who verified \eqref{eqn:dlvp} for $R\approx 5.573$ and $T=2$.

\textit{Koborov--Vinogradov zero-free region.}
In 1958, Koborov \cite{koborov58} and Vinogradov \cite{vinogradov58} independently demonstrated that there exists a positive constant $R_1$ such that $\zeta$ is non-zero in the region $s=\sigma + it$ such that $t\geq T$ and
\begin{equation}\label{eqn:kob-vino}
    \sigma \geq 1-\frac{1}{R_1(\log t)^{\frac{2}{3}}(\log\log t)^{\frac{1}{3}}}.
\end{equation}
The best known zero-free region for $\zeta$ of this kind is attributed to Ford \cite{ford2002zero}, who has verified \eqref{eqn:kob-vino} for $R_1 = 57.54$ and $T=3$. Ford \cite{ford2002zero} also establishes the zero-free region \eqref{eqn:kob-vino} for large $t$ with $R_1 = 49.13$.

Naturally, the closest form of the zero-free region for $\zeta_L$ will also depend on the extra variables $d_L$ and $n_L$. However, the method we adopt is based on de la Vall{\'e}e Poussin's method for determining the classical zero-free region for $\zeta$. One complication is that a so-called \textit{exceptional} zero could exist inside a zero-free region for $\zeta_L$. If this exceptional zero exists, then {it} must be simple and real.

Kadiri \cite[Theorem 1.1]{kadiri2012explicit} was the last to re-purpose de la Vall\'{e}e Poussin's proof (using Ste\v{c}kin's \cite{Stechkin1970} so-called differencing trick) to obtain a zero-free region for $\zeta_L$. In this paper, we will establish Theorem \ref{thm:main_result1}, a new zero-free region for $\zeta_L$ which builds upon Kadiri's zero-free region for $\zeta_L$. We will also establish Theorem \ref{thm:main_result2}, which will reveal more information pertaining to the exceptional zero.

\begin{thm}\label{thm:main_result1}
Suppose $(C_1,C_2,C_3,C_4) = (12.2411, 9.5347, 0.05017, 2.2692)$, then $\zeta_L(\sigma + it)$ is non-zero for
\begin{equation}\label{eqn:main_result1}
    \sigma \geq 1 - \frac{1}{C_1\log d_L + C_2 \cdot n_L\log |t|+ C_3\cdot n_L + C_4} \mbox{ and }|t|\geq 1.
\end{equation}
\end{thm}

\begin{thm}\label{thm:main_result2}
For asymptotically large $d_L$ and $R = 12.43436$, $\zeta_L(\sigma + it)$ has at most one zero in the region
\begin{equation}\label{eqn:main_result2}
    \sigma \geq 1 - \frac{1}{R \log d_L}\mbox{ and }|t| < 1
\end{equation}
If this exceptional zero exists, then it is simple and real.
\end{thm}

Kadiri \cite{kadiri2012explicit} established \eqref{eqn:main_result1} with $(C_1, C_2, C_3, C_4) = (12.55, 9.69, 3.03, 58.63)$.
To yield Theorem \ref{thm:main_result1}, we will follow a similar process to Kadiri, but observe some improvements. An important step in the proof of Theorem \ref{thm:main_result1} is to choose a polynomial $p_n(\varphi)$ from the so-called the class of non-negative, trigonometric polynomials of degree $n$; denoted and defined by
$$P_n:=\left\{p_n(\varphi)=\sum_{k=0}^n a_k\cos (k\varphi):p_n(\varphi)\geq 0\mbox{ for all }\varphi\mbox{, }a_k\geq 0\mbox{ and }a_0<a_1\right\}.$$
Whereas Kadiri worked with polynomials from $P_4$, we will use the same polynomial from $P_{16}$ as Mossinghoff--Trudgian \cite{MossinghoffTrudgian2015}. This polynomial has been optimised by simulated annealing for computations pertaining to their computations for the zero-free region for $\zeta$. This amendment contributed \textit{all} of the improvements that can be seen for $C_1$ and $C_2$. In fact, if one re-runs Kadiri's computations, only updating the polynomial, then this establishes \eqref{eqn:main_result1} with $(C_1, C_2, C_3, C_4) = (12.2411, 9.5347, 3.3492, 57.7027)$.

Another improvement follows from improvements we have made to \cite[Lemma 2]{mccurley84} from McCurley. In particular, we improve explicit values for $\mcS(k)$, a computable constant dependent on $k\in\NN$. These improvements will contribute almost all of the improvement one observes for $C_3$.

Kadiri \cite{kadiri2012explicit} also established \eqref{eqn:main_result2} with $R = 12.7305$. To yield Theorem \ref{thm:main_result2}, we will recycle bounds from \cite[§3]{kadiri2012explicit} and apply the same higher degree polynomial from $P_{16}$.
A corollary of the method we use to establish Theorem \ref{thm:main_result2} is an improvement to a well-known region by Stark \cite{stark1974some}. However, because we only update the polynomial for this method, we cannot improve Stark's result further than \cite[Corollary 1.2]{kadiri2012explicit} already does.

Finally, if an exceptional zero $\beta_1$ exists, then one can enlarge the zero-free region in Theorem \ref{thm:main_result2} using the Deuring-Heilbronn phenomenon \cite{linnik_44_2}. This was one of the key ingredients in work by Ahn--Kwon \cite{ahn2019explicit}, Zaman \cite{zaman2017bounding} and Kadiri--Ng--Wong \cite{kadiri_ng_wong2019}, which pertains to the least prime ideal in the Chebotar\"{e}v density theorem.

\textit{Remark.}
The method of proof which we follow does not use Heath-Brown’s version of Jensen’s formula \cite[Lemma 3.2]{heath1992zero}, although this might yield better zero-free regions than those we can obtain using this method. This is partially because there does \textit{not} exist a \textit{general} sub-convexity bound for \textit{general} number fields, so it is difficult to apply his approach in the number field setting --- see Kadiri \cite{kadiri2012explicit} for an excellent explanation of this.

\subsection*{Acknowledgements}

I would like to thank my supervisor, Tim Trudgian, for bringing this project to my attention, as well as his continued support. I would also like to thank my other colleagues at UNSW Canberra for their support throughout this process and Kevin Ford for correcting a referential error.
\section{Proof of Theorem \ref{thm:main_result1}}

The set-up of our proof for Theorem \ref{thm:main_result1} is the same as that which Kadiri uses in her proof of \cite[Theorem 1.1]{kadiri2012explicit}, which has a similar shape to Ste\v{c}kin's argument \cite{Stechkin1970} for $\zeta$. Suppose $t \geq 1$. We introduce some definitions, which will hold for the remainder of this paper:
\begin{itemize}
    \item $\kappa = \frac{1}{\sqrt{5}}$;
    \item $s_k = \sigma + ikt$ such that $k\in \NN$, $1<\sigma < 1 + \varepsilon$ for some $0 < \varepsilon \leq 0.15$;
    \item $s'_k = \sigma_1 + ikt$ such that $\sigma_1 = \frac{1 + \sqrt{1 + 4\sigma^2}}{2}$.
\end{itemize}
{Note that $\sigma_1$ depends on $\sigma$, so for convenience we will write $\sigma_1(a)$ to denote the the value of $\sigma_1$ at $\sigma = a$.}
{To prove Theorem \ref{thm:main_result1}, we will isolate a non-trivial zero $\rho = \beta + it$ of $\zeta_L$ such that $\beta > 1 - \varepsilon \geq 0.85$, choose a polynomial $p_{n}(\varphi)$ from $P_{n}$, and} consider the function
$$S(\sigma,t) = \sum_{k=0}^{n} a_kf_L(\sigma, kt),$$ such that
\begin{align*}
    f_L(\sigma,kt)
    &= -\Rea\left(\frac{\zeta'_L}{\zeta_L}(s_k) - \kappa \frac{\zeta'_L}{\zeta_L}({s'_k})\right)\\
    &= \sum_{0\neq\mfP\subset\mcO_L} \Lambda(\mfP)(N(\mfP)^{-\sigma}-\kappa N(\mfP)^{-\sigma_1}) \cos(k t\log(N(\mfP)).\label{eqn:useful1}
\end{align*}
It follows that
$$S(\sigma,t)= \sum_{0\neq\mfP\subset\mcO_L} \Lambda(\mfP)(N(\mfP)^{-\sigma}-\kappa N(\mfP)^{-\sigma_1}) p_{n}(t\log(N(\mfP))\geq 0.$$
On the other hand, we can utilise the explicit formula \cite[(8.3)]{lagarias1977effective},
\begin{equation}\label{eqn:explicit_formula}
    - \frac{\zeta_L'}{\zeta_L}(s_k) = \frac{\log d_L}{2} + \frac{1}{s_k} + \frac{1}{s_k - 1} + \frac{\gamma_L'}{\gamma_L}(s_k) - \frac{1}{2}\sum_{\varrho\in Z(\zeta_L)} \left(\frac{1}{s_k - \varrho} + \frac{1}{s_k - \overline{\varrho}}\right).
\end{equation}
Here, $Z(\zeta_L)$ denotes the set of non-trivial zeros of $\zeta_L$. One can use \eqref{eqn:explicit_formula} to show
\begin{equation}\label{eqn:imp}
    0\leq S(\sigma,t)\leq S_1+S_2+S_3+S_4,
\end{equation}
where $F(s,z) = \Rea\left(\frac{1}{s-z}+\frac{1}{s-1+\bar{z}}\right)$ such that
\begin{align*}
    S_1 &= -\sum_{k=0}^{n} a_k\sum_{\varrho\in Z(\zeta_L)}\Rea\left(\frac{1}{s_k-\varrho} - \frac{\kappa}{{s'_k}-\varrho}\right),\\
    S_2 &= \frac{1-\kappa}{2} \left(\sum_{k=0}^{n} a_k\right) \log d_L,\\
    S_3 &= \sum_{k=0}^{n} a_k \left(F(s_k,1) -\kappa F({s'_k},1)\right),\mbox{ and}\\
    S_4 &= \sum_{k=0}^{n} a_k \Rea\left(\frac{\gamma_L'(s_k)}{\gamma_L(s_k)} - \kappa \frac{\gamma_L'({s'_k})}{\gamma_L({s'_k})}\right).
\end{align*}
{We will choose $n=16$, so that we can apply Mossinghoff--Trudgian's polynomial $p_{16}(\varphi)\in P_{16}$ from \cite{MossinghoffTrudgian2015}. Taking $n=16$, $S_2$ is directly computable, and we find upper bounds for $S_1$, $S_3$, and $S_4$ in Sections \ref{subsec:S1}, \ref{subsec:S3}, and \ref{subsec:S4}. The resulting upper bound for $S_1 + S_2 + S_3 + S_4$ will depend on $\beta$, $\sigma$, $t$, the coefficients of $p_{16}(\varphi)$ and $\varepsilon$, therefore we may use \eqref{eqn:imp} and rearrange the inequality to obtain Theorem \ref{thm:main_result1} in Section \ref{subsec:computations}.}

\subsection{Upper bound for $S_1$}\label{subsec:S1}

\begin{lemma}[Ste\v{c}kin \cite{Stechkin1970}]\label{lem:Stechkin_lemma}
Suppose $s=\sigma + it$ with $1 < \sigma \leq 1.25$ and $z\in\CC$. If $0 < \Rea(z) < 1$, then
\begin{equation}\label{Stechkin_lemma:ineq_1}
    F(s, z)-\kappa F({s'_1},z)\geq 0.
\end{equation}
Moreover, if $\Ima(z)=\Ima(s)=t$ and $\frac{1}{2} \leq \Rea(z) < 1$, then
\begin{equation*}
    \Rea\left(\frac{1}{s-1+\bar{z}}\right)-\kappa F({s'_1},z)\geq 0.
\end{equation*}
\end{lemma}

Note that $\kappa$ is the largest value such that \eqref{Stechkin_lemma:ineq_1} holds. This subsection is \textit{not} an improvement on \cite[Lemma 2.3]{kadiri2012explicit}, rather a repeat for the purpose of clarity. By the positivity condition \eqref{Stechkin_lemma:ineq_1} in Lemma \ref{lem:Stechkin_lemma}, we have
\begin{equation}\label{eqn:ell_ineq}
    \ell(s_k) := \sum_{\varrho\in Z(\zeta_L)}\Rea\left(\frac{1}{s_k-\varrho} - \frac{\kappa}{{s'_k}-\varrho}\right) \leq \kappa F({s'_k}, \rho) - F(s_k,\rho).
\end{equation}
If $k = 1$, then \eqref{eqn:ell_ineq} implies that
$$\ell(s_1)\leq - \frac{1}{\sigma - \beta} - \frac{1}{\sigma - 1 + \beta} + \frac{\kappa}{\sigma_1 - \beta} + \frac{\kappa}{\sigma_1 - 1 + \beta} = - \frac{1}{\sigma - \beta} + g(\sigma, \beta).$$
We see that $g(\sigma, \beta) < g(1, 1)$ and $g(1, 1)$ is small and negative, so $\ell(s_1) \leq - \frac{1}{\sigma - \beta}$. Moreover, if $k \neq 1$, then \eqref{eqn:ell_ineq} implies that $\ell(s_k)\leq 0$ by \eqref{Stechkin_lemma:ineq_1}. One can package the preceding observations into the following lemma.

\begin{lemma}\label{lem:S1}
Isolate a zero $\rho = \beta + it \in Z(\zeta_L)$ such that $\beta \geq 1 - \varepsilon \geq 0.85$, then
$$\ell(\sigma + ikt) \leq
\begin{cases}
- \frac{1}{\sigma - \beta} &\mbox{if }k = 1,\\
0&\mbox{if }k\neq 1.
\end{cases}$$
Therefore, $S_1 \leq - \frac{a_1}{\sigma - \beta}.$
\end{lemma}

\subsection{Upper bound for $S_3$}\label{subsec:S3}

Suppose that
\begin{align*}
    \Sigma_k(\sigma,t)
    &:= F(\sigma + ikt, 1) - \kappa F(\sigma_1 + ikt, 1)\\
    &= \frac{\sigma}{\sigma^2 + k^2 t^2} + \frac{\sigma - 1}{(\sigma - 1)^2 + k^2 t^2} - \kappa\frac{\sigma_1}{{\sigma_1}^2 + k^2 t^2} - \kappa \frac{\sigma_1 - 1}{(\sigma_1 - 1)^2 + k^2 t^2}.
\end{align*}

\textbf{Case I.}
If $k = 0$, then $\Sigma_k$ is only dependent on $\sigma$, with a singularity occuring at $\sigma = 1$. In fact,
$$\Sigma_0(\sigma,t) = \frac{1}{\sigma} + \frac{1}{\sigma - 1} - \frac{\kappa}{\sigma_1} - \frac{\kappa}{\sigma_1 - 1} :=  \frac{1}{\sigma - 1} + h(\sigma).$$
We observe that $h(\sigma)$ \textit{increases} as $\sigma$ increases, so for $\alpha_{\varepsilon} = h(1 + \varepsilon) < 0.021467$, we have
$$\Sigma_0(\sigma,t) \leq \frac{1}{\sigma - 1} + \alpha_{\varepsilon}.$$

\textbf{Case II.}
Suppose $1 \leq k \leq 16$, then $\Sigma_k(\sigma,t)$ depends on $\sigma$ and $t$. For each $\sigma$, $\Sigma_k(\sigma,t)$ \textit{decreases} as $t$ increases, {because the derivative of $\Sigma_k(\sigma,t)$ with respect to $t$ is negative for all $t\geq 1$}.
Therefore, $\Sigma_k(\sigma,t)\leq \Sigma_k(\sigma,1)$, which \textit{increases} as $\sigma$ increases, {because the derivative of $\Sigma_k(\sigma,1)$ with respect to $\sigma$ is positive for all $1\leq\sigma\leq 1.15$}.
It follows that
$$\Sigma_k(\sigma,t) \leq \Sigma_k(1+\varepsilon,1) < \mathcal{B}_{\varepsilon}(k),$$
where admissible values for $\mathcal{B}_{\varepsilon}(k)$ are easily computed using a computer. {To further verify this bound, the Maximize command in Maple confirms that the maximum of $\Sigma_k(\sigma,t)$ occurs at $\sigma = 1+\varepsilon$ and $t=1$.} For example, if $\varepsilon = 0.15$ or $\varepsilon = 0.01$, then admissible values of $\mathcal{B}_{0.15}(k)$ and $\mathcal{B}_{0.01}(k)$ are given in Table \ref{tab:admissible_valus_B_k_0.15} and Table \ref{tab:admissible_valus_B_k_0.01} respectively. Note that we round up at 8 decimal places, to account for any possible rounding errors.

\begin{table}[!htb]
\begin{minipage}{.5\linewidth}
    \centering
    \begin{tabular}{cccc}
    $k$ & $\mathcal{B}_{0.15}(k)$ & $k$ & $\mathcal{B}_{0.15}(k)$\\
    \hline
    1 & 0.23445352 & 9 & 0.00235718 \\
    2 & 0.06869804 & 10 & 0.00188669 \\
    3 & 0.02783858 & 11 & 0.00154513 \\
    4 & 0.01427867 & 12 & 0.00128917 \\
    5 & 0.0085573  & 13 & 0.0010924 \\
    6 & 0.00568194 & 14 & 0.00093759 \\
    7 & 0.00404715 & 15 & 0.00081374 \\
    8 & 0.00303134 &  16 & 0.00071303
    \end{tabular}
    \caption{Admissible values for $\mathcal{B}_{0.15}(k)$.}
    \label{tab:admissible_valus_B_k_0.15}
\end{minipage}%
\begin{minipage}{.5\linewidth}
    \centering
    \begin{tabular}{cccc}
    $k$ & $\mathcal{B}_{0.01}(k)$ & $k$ & $\mathcal{B}_{0.01}(k)$\\
    \hline
    1 & 0.10919579  & 9 & 0.00029396\\
    2 & 0.03040152  & 10& 0.00021655\\
    3 & 0.00958566 & 11 & 0.00016557\\
    4 & 0.00384196 & 12 & 0.00013046\\
    5 & 0.00185609 & 13 & 0.00010535\\
    6 & 0.00102853 & 14 & 0.00008684\\
    7 & 0.00063099 & 15 & 0.00007282\\
    8 & 0.00041809 & 16 & 0.00006196
    \end{tabular}
    \caption{Admissible values for $\mathcal{B}_{0.01}(k)$.}
    \label{tab:admissible_valus_B_k_0.01}
\end{minipage} 
\end{table}

Now, we can collect the preceding observations to yield Lemma \ref{lem:S3}.

\begin{lemma}\label{lem:S3}
For $0\leq k\leq 16$, we have that
$$\Sigma_k(\sigma,t) \leq
\begin{cases}
\frac{1}{\sigma - 1} + \alpha_{\varepsilon} &\mbox{if }k = 0,\\
\mathcal{B}_{\varepsilon}(k) &\mbox{if }k\neq 0.
\end{cases}$$
Under a choice of polynomial from $P_{16}$, it follows that
$$S_3 \leq a_0\left(\frac{1}{\sigma - 1} + \alpha_{\varepsilon}\right) + \sum_{k = 1}^{16} a_k \mathcal{B}_{\varepsilon}(k).$$
\end{lemma}

\textit{Remark.}
The benefits of Lemma \ref{lem:S3} over \cite[Lemma 2.4]{kadiri2012explicit} lie in the computed constants $\mathcal{B}_{\varepsilon}(k)$. That is, Kadiri established $\Sigma_k(\sigma,t)\leq 1.6666$ for $1\leq k\leq 4$.

\subsection{Upper bound for $S_4$}\label{subsec:S4}

We bring forward an observation from Kadiri \cite[§2.4]{kadiri2012explicit},
\begin{align}
    &\Rea\left(\frac{\gamma_L'(s_k)}{\gamma_L(s_k)} - \kappa \frac{\gamma_L'({s'_k})}{\gamma_L({s'_k})}\right)
    \leq -\frac{1-\kappa}{2}\cdot \log\pi \cdot n_L \nonumber\\
    &\quad\qquad\qquad +\frac{n_L}{2}\max_{\delta\in\{0,1\}}\left\{\Rea\left(\frac{\Gamma'}{\Gamma}\left(\frac{s_k + \delta}{2}\right) - \kappa \frac{\Gamma'}{\Gamma}\left(\frac{{s'_k} + \delta}{2}\right)\right)\right\}.\nonumber
\end{align}

\textbf{Case I.}
If $k = 0$, then we directly compute that
\begin{equation}\label{eqn:justifications1}
    \frac{1}{2}\max_{\delta\in\{0,1\}}\left\{\Rea\left(\frac{\Gamma'}{\Gamma}\left(\frac{\sigma + \delta}{2}\right) - \kappa \frac{\Gamma'}{\Gamma}\left(\frac{\sigma_1 + \delta}{2}\right)\right)\right\} \leq d_{\varepsilon}(0),
\end{equation}
where $d_{\varepsilon}(0)$ is the maximum of the functions such that $\sigma = 1 + \varepsilon$. To see this, one can observe that the left-hand side of \eqref{eqn:justifications1} is maximised at $\sigma = 1 + \varepsilon$ visually \textit{or} use the Maximize command in Maple. For example, if $\varepsilon = 0.01$, then
$$d_{0.01}(0) = -0.2500763736.$$

\textbf{Case II.}
Suppose $1 \leq k\leq 16$. McCurley  \cite[Lemma 2]{mccurley84} establishes that
\begin{align}
    &\frac{1}{2}\Rea\left(\frac{\Gamma'}{\Gamma}\left(\frac{s_k + \delta}{2}\right) - \kappa \frac{\Gamma'}{\Gamma}\left(\frac{{s'_k} + \delta}{2}\right)\right)
    = \frac{1 - \kappa}{2}\log\frac{kt}{2}
    + \Xi(\sigma,k,t,\delta)\nonumber\\
    &\hspace{0.5cm}+ \frac{\theta_1}{2k} \left(\frac{\pi}{2} - \arctan\left(\frac{1 + \delta}{k}\right)\right)
    + \kappa\frac{\theta_2}{2k} \left(\frac{\pi}{2} - \arctan\left(\frac{\sigma_1(1) + \delta}{k}\right)\right),\label{eqn:McCurley}
\end{align}
where $|\theta_i|\leq 1$ and
\begin{align}
    \Xi(\sigma,k,t,\delta)
    &= \frac{1}{4}\log\left[1 + \left(\frac{\sigma + \delta}{kt}\right)^2\right]
    - \frac{\kappa}{4}\log\left[1 + \left(\frac{\sigma_1 + \delta}{kt}\right)^2\right]\nonumber\\
    &\hspace{1.5cm}- \frac{\sigma + \delta}{2((\sigma + \delta)^2 + k^2t^2)}
    +\kappa \frac{\sigma_1 + \delta}{2((\sigma_1 + \delta)^2 + k^2t^2)}.\nonumber
\end{align}
Next, we will bound $\Xi(\sigma,k,t,\delta)$ using two different methods, then choose the best bound for each $k$.

\textit{Method I.}
For any $t > 0$, we have
\begin{align*}
    \Xi_1(\sigma,k,t,\delta)
    &:= - \frac{\sigma + \delta}{2((\sigma + \delta)^2 + k^2t^2)} +\kappa \frac{\sigma_1 + \delta}{2((\sigma_1 + \delta)^2 + k^2t^2)}\\
    &\leq \frac{\kappa (\sigma_1 + \delta) - \sigma - \delta}{2((\sigma_1 + \delta)^2 + k^2t^2)}
    \leq 0,
\end{align*}
because $\sigma < \sigma_1$ and $\kappa (\sigma_1 + \delta) - \sigma - \delta \leq 0$.
Moreover, for fixed $\sigma$, observe that
\begin{align*}
    \Xi_2(\sigma,k,t,\delta)
    := \frac{1}{4}\log\left[1 + \left(\frac{\sigma + \delta}{kt}\right)^2\right] - \frac{\kappa}{4}\log\left[1 + \left(\frac{\sigma_1 + \delta}{kt}\right)^2\right]
\end{align*}
is positive for $t\geq 1$ and decreases as $t$ increases, because the derivative of $\Xi_2(\sigma,k,t,\delta)$ with respect to $t$ is negative for all $t\geq 1$. Therefore,
\begin{equation*}
    \Xi_2(\sigma,k,t,\delta)
    \leq \Xi_2(\sigma,k,1,\delta)
\end{equation*}
for $t\geq 1$, which increases as $\sigma$ increases in the range $1 \leq \sigma \leq 1.15$, because the derivative of $\Xi_2(\sigma,k,1,\delta)$ with respect to $\sigma$ is positive for $1 \leq \sigma \leq 1.15$. Hence,  for each $k$,
\begin{equation*}
    \Xi_2(\sigma,k,t,\delta)
    \leq \Xi_2(1+\varepsilon,k,1,\delta).
\end{equation*}
To verify the preceding bound, the Maximize command in Maple confirms that the maximum of $\Xi_2(\sigma,k,t,\delta)$ occurs at $\sigma = 1+\varepsilon$ and $t=1$.
It follows that $\Xi(\sigma,k,t,\delta) \leq \Xi_2(1+\varepsilon,k,1,\delta)$ for each $k$ and
$$\frac{1}{2}\max_{\delta\in\{0,1\}}\left\{\Rea\left(\frac{\Gamma'}{\Gamma}\left(\frac{s_k + \delta}{2}\right) - \kappa \frac{\Gamma'}{\Gamma}\left(\frac{{s'_k} + \delta}{2}\right)\right)\right\}
\leq \frac{1 - \kappa}{2}\log t + \mcS_1(k,\varepsilon),$$
where $\mcS_1(k,\varepsilon) = \max_{\delta\in\{0,1\}}\left\{\mcC_1(k,\delta,\varepsilon)\right\}$ such that
\begin{align*}
    \mcC_1(k,\delta,\varepsilon)
    := &\frac{1 - \kappa}{2}\log \frac{k}{2} + \Xi_2(1+\varepsilon,k,1,\delta)\\
    &+ \frac{1}{2k} \left(\frac{\pi}{2} - \arctan\left(\frac{1 + \delta}{k}\right)\right) + \frac{\kappa}{2k} \left(\frac{\pi}{2} - \arctan\left(\frac{\sigma_1(1) + \delta}{k}\right)\right).
\end{align*}

\textit{Method II.}
We will verify that for $0<\varepsilon\leq 0.15$,
\begin{equation}\label{eqn:revison_justification}
    \Xi(\sigma,k,t,\delta)\leq \mcA(k,\delta,\varepsilon) :=
    \begin{cases}
        0&\mbox{if }\delta = 0,\\
        0&\mbox{if }\delta = 1\mbox{ and }k\not\in\{1,2\},\\
        \Xi(1+\varepsilon,k,1,1)&\mbox{if }\delta = 1\mbox{ and }k = 1.\\
        \Xi(1.15,k,1,1)&\mbox{if }\delta = 1\mbox{ and }k=2.
\end{cases}
\end{equation}
First, for fixed $\sigma$ and $\delta = 0$, the derivative of $\Xi(\sigma,k,t,\delta)$ with respect to $t$ is positive for $t\geq 1$, so $\Xi(\sigma,k,t,0)$ is increasing as $t\to\infty$. Therefore, for each $\sigma\in [1,1.15]$,
\begin{align*}
    \Xi(\sigma,k,t,0) \leq \lim_{t\to\infty}\Xi(\sigma,k,t,0) = 0.
\end{align*}
Next, for fixed $\sigma$ and $\delta = 1$, the derivative of $\Xi(\sigma,k,t,\delta)$ with respect to $t$ is positive for $t\geq 1$ whenever $k\not\in \{1,2,3\}$, so $\Xi(\sigma,k,t,1)$ is increasing as $t\to\infty$ for $k\not\in \{1,2,3\}$. Therefore, for each $k\not\in \{1,2,3\}$ and $1\leq \sigma \leq 1.15$,
\begin{align*}
    \Xi(\sigma,k,t,1) \leq \lim_{t\to\infty}\Xi(\sigma,k,t,1) = 0.
\end{align*}
To completely verify \eqref{eqn:revison_justification}, we now establish bounds for the special cases $\delta = 1$ and $k\in \{1,2,3\}$. Observe that for each $t\geq 1$, the derivative of $\Xi(\sigma,k,t,1)$ with respect to $\sigma$ is positive for $1\leq \sigma \leq 1.15$ whenever $k\in \{1,2,3\}$, so
\begin{equation}\label{eqn:in_bed0}
    \Xi(\sigma,k,t,1) \leq \Xi(1+\varepsilon,k,t,1).
\end{equation}
Suppose that $k\in\{1,2,3\}$ and observe that in the range $t\geq 1$, $\Xi(1+\varepsilon,k,t,1)$ \textit{either} has one minimum point at $t=t_k(\varepsilon)$ \textit{or} increases as $t\to\infty$. Here, $t_k(\varepsilon)$ equals the only root of the derivative of $\Xi(1+\varepsilon,k,t,1)$ with respect to $t$ in the range $t\geq 1$. If this root does not exist, then set $t_k(\varepsilon) = 1$ for convenience. For example, $t_1(0.15) =3.2308\ldots$, $t_2(0.15)=1.6154\ldots$, $t_3(0.15) = 1.0769\ldots$ and $t_3(0.01) = 1$. It follows that $\Xi(1+\varepsilon,k,t,1)$ decreases for $1\leq t \leq t_k(\varepsilon)$ and $\Xi(1+\varepsilon,k,t,1)$ increases for $t > t_k(\varepsilon)$, so
\begin{align*}
    \Xi(1+\varepsilon,k,t,1)
    \leq
    \begin{cases}
    \Xi(1+\varepsilon,k,1,1) & \text{if }1\leq t\leq t_k(\varepsilon),\\
    \lim_{t\to\infty}\Xi(1+\varepsilon,k,t,1) & \text{if }t> t_k(\varepsilon),
    \end{cases}
\end{align*}
in which $\lim_{t\to\infty}\Xi(1+\varepsilon,k,t,1) = 0$ for each $k$.
If $k=1$, then for $t\geq 1$, we have
\begin{align}
    \Xi(1+\varepsilon,1,t,1)
    \leq \max\left\{\Xi(1+\varepsilon,1,1,1),0\right\}
    = \Xi(1+\varepsilon,1,1,1).\label{eqn:in_bed1}
\end{align}
Observe that $\Xi(1+\varepsilon,2,1,1)$ increases as $0<\varepsilon\leq 0.15$ increases. So, if $k=2$, then for $t\geq 1$, we have
\begin{align}
    \Xi(1+\varepsilon,2,t,1)
    \leq \max\left\{\Xi(1+\varepsilon,2,1,1),0\right\}
    &\leq \max\left\{\Xi(1.15,2,1,1),0\right\}\nonumber\\
    &\leq \Xi(1.15,2,1,1).\label{eqn:in_bed2}
\end{align}
In this case, the final bound is convenient and not too wasteful, because $\Xi(1.15,2,1,1)$ is small. Finally, if $k=3$, then for $t\geq 1$, we have
\begin{equation}\label{eqn:in_bed3}
    \Xi(1+\varepsilon,3,t,1)
    \leq \max\left\{\Xi(1+\varepsilon,3,1,1),0\right\}
    = 0.
\end{equation}
Combining the observation \eqref{eqn:in_bed0} with \eqref{eqn:in_bed1}, \eqref{eqn:in_bed2}, and \eqref{eqn:in_bed3} implies \eqref{eqn:revison_justification}.
For each $k$, it follows from \eqref{eqn:McCurley} and \eqref{eqn:revison_justification} that
$$\frac{1}{2}\max_{\delta\in\{0,1\}}\left\{\Rea\left(\frac{\Gamma'}{\Gamma}\left(\frac{s_k + \delta}{2}\right) - \kappa \frac{\Gamma'}{\Gamma}\left(\frac{{s'_k} + \delta}{2}\right)\right)\right\}
\leq \frac{1 - \kappa}{2}\log t + \mcS_2(k,\varepsilon),$$
where $\mcS_2(k,\varepsilon) = \max_{\delta\in\{0,1\}}\left\{\mcC_2(k,\delta,\varepsilon)\right\}$ such that
\begin{align*}
    \mcC_2(k,\delta,\varepsilon)
    := &\frac{1 - \kappa}{2}\log \frac{k}{2} + \mcA(k,\delta,\varepsilon)\\
    &+ \frac{1}{2k} \left(\frac{\pi}{2} - \arctan\left(\frac{1 + \delta}{k}\right)\right) + \frac{\kappa}{2k} \left(\frac{\pi}{2} - \arctan\left(\frac{\sigma_1(1) + \delta}{k}\right)\right).
\end{align*}

\textit{Combination.}
We say that $\mcS(k,\varepsilon)=\min(\mcS_1(k,\varepsilon), \mcS_2(k,\varepsilon))$ and (for $1\leq k \leq 16$) present the quantities $\mcS_1(k,0.15)$, $\mcS_2(k,0.15)$ and $\mcS(k,0.15)$ alongside each other in Table \ref{tab:values_for_S4}. It turns out that $\mcS_2(k,0.15)$ yields a better bound for cases $k=1,2,3,4$ and $\mcS_1(k,0.15)$ yields the better bound otherwise. Finally, we package our observations into a useful lemma (Lemma \ref{lem:S4}).

\begin{table}[h!]
    \centering
    \begin{tabular}{c|c|c|c}
         $k$ & $\mcS_1(k,0.15)$ & $\mcS_2(k,0.15)$ & $\mcS(k,0.15)$ \\
         \hline
          1& 0.3784516540 & 0.3249009026 & 0.3249009026\\
          2& 0.3839873212 & 0.3763572015 & 0.3763572015\\
          3& 0.4018562060 & 0.4004551145 & 0.4004551145\\
          4& 0.4238223974 & 0.4236306767 & 0.4236306767\\
          5& 0.4467597648 & 0.4468482525 & 0.4467597648\\
          6& 0.4693610537 & 0.4695098183 & 0.4693610537\\
          7& 0.4910902618 & 0.4912403488 & 0.4910902618\\
          8& 0.5117562107 & 0.5118920810 & 0.5117562107\\
          9& 0.5313238925 & 0.5314428586 & 0.5313238925\\
          10& 0.5498280118 & 0.5499312088 & 0.5498280118\\
          11& 0.5673323540 & 0.5674218683 & 0.5673323540\\
          12& 0.5839104248 & 0.5839883668 & 0.5839104248\\
          13& 0.5996362678 & 0.5997044990 & 0.5996362678\\
          14& 0.6145802698 & 0.6146403531 & 0.6145802698\\
          15& 0.6288074426 & 0.6288606647 & 0.6288074426\\
          16& 0.6423769295 & 0.6424243440 & 0.6423769295

    \end{tabular}
    \caption{Computed values for $\mcS_1(k,0.15)$, $\mcS_2(k,0.15)$ and $\mcS(k,0.15)$.}
    \label{tab:values_for_S4}
\end{table}

\begin{lemma}\label{lem:S4}
For $0\leq k \leq 16$, we have shown that
$$\Rea\left(\frac{\gamma_L'(s_k)}{\gamma_L(s_k)} - \kappa \frac{\gamma_L'({s'_k})}{\gamma_L({s'_k})}\right) \leq 
\begin{cases}
n_L\left(d_{\varepsilon}(0) -\frac{1-\kappa}{2}\cdot \log\pi \right)&\mbox{if }k=0,\\
n_L\left(\frac{1 - \kappa}{2}\left(\log t + \log\left(\frac{k}{\pi}\right)\right) + \mcS(k,\varepsilon) \right)&\mbox{if }k\neq 0.
\end{cases}$$
Under a choice of polynomial from $P_{16}$, it follows that
\begin{align*}
    S_4 
    \leq a_0 n_L&\left(d_{\varepsilon}(0) -\frac{1-\kappa}{2}\cdot \log\pi \right)\\
    &\hspace{1.75cm}+ \sum_{k = 1}^{16} a_k n_L\left(\frac{1 - \kappa}{2}\left(\log t + \log\left(\frac{k}{\pi}\right)\right) + \mcS(k,\varepsilon) \right).
\end{align*}
\end{lemma}

%\textit{Remark.}
%Using the same methods, we expect that one can extend Lemma \ref{lem:S4} to hold for $k > 16$ without any trouble.

\textit{Remark.}
The benefits of Lemma \ref{lem:S4} over \cite[Lemma 2.5]{kadiri2012explicit} lie in the computed constants $d_\varepsilon(0)$ and $\mathcal{S}(k,\varepsilon)$. Kadiri imports results from McCurley \cite[Lemma 2]{mccurley84} for her bound, so the improvements we see follow from our observations pertaining to McCurley's work.

\subsection{Computations}\label{subsec:computations}

\begin{table}[h]
\centering
\begin{tabular}{c|l}
     $a_0$  &  $1$\\
     $a_1$  &  $1.74126664022806$\\
     $a_2$  &  $1.128282822804652$\\
     $a_3$  &  $0.5065272432186642$\\
     $a_4$  &  $0.1253566902628852$\\
     $a_5$  &  $2.372710620\cdot 10^{-26}$\\
     $a_6$  &  $2.818732841\cdot 10^{-22}$\\
     $a_7$  &  $0.01201214561729989$\\
     $a_8$  &  $0.006875849760911001$\\
     $a_9$  &  $2.064157910\cdot 10^{-23}$\\
     $a_{10}$  &  $6.601587090\cdot 10^{-11}$\\
     $a_{11}$  &  $0.001608306592372963$\\
     $a_{12}$  &  $0.001017994683287104$\\
     $a_{13}$  &  $6.728831293\cdot 10^{-11}$\\
     $a_{14}$  &  $3.682448595\cdot 10^{-11}$\\
     $a_{15}$  &  $2.949853019\cdot 10^{-6}$\\
     $a_{16}$  &  0.00003713656497\\
\end{tabular}
\caption{Table of coefficients for Mossinghoff--Trudgian's polynomial $p_{16}(\varphi)\in P_{16}$.}
\label{table:cefficients1}
\end{table}

As declared in the introduction, we will choose the polynomial $p_{16}(\varphi)\in P_{16}$ from \cite{MossinghoffTrudgian2015}, whose coefficients are given in Table \ref{table:cefficients1}. Suppose $r>0$ and $\sigma$ is chosen such that $\sigma - 1 = r(1-\beta)$ where $\rho = \beta + it \in Z(\zeta_L)$ is an isolated zero such that $\beta \geq 1 - \varepsilon \geq 0.85$. Applying the upper bounds for each $S_i$, which can be found in {Lemmas} \ref{lem:S1}, \ref{lem:S3} and \ref{lem:S4}, then rearranging inequality \eqref{eqn:imp} will yield 
\begin{equation}\label{eqn:useful2}
    \beta\leq 1 - \frac{\frac{a_1}{1+r} - \frac{a_0}{r}}{c_1\log d_L + c_2 n_L \log t + c_3 n_L + c_4},
\end{equation}
where
\begin{align*}
    c_1 &= \frac{1-\kappa}{2}\sum_{k=0}^{16} a_k,\\
    c_2 &= \frac{1-\kappa}{2}\sum_{k=1}^{16} a_k,\\
    c_3 &= a_0 \left(d_{\varepsilon}(0) - \frac{1-\kappa}{2}\log \pi\right)+ \sum_{k=1}^{16} a_k\left(\frac{1 - \kappa}{2}\log\left(\frac{k}{\pi}\right) + \mcS(k,\varepsilon) \right)\mbox{ and}\\
    c_4 &= \alpha_{\varepsilon} a_0 + \sum_{k=1}^{16} a_k \mathcal{B}_{\varepsilon}(k).
\end{align*}
For the remainder of this proof, we replicate the process which Kadiri \cite{kadiri2012explicit} followed. The maximum value of $\frac{a_1}{1+r} - \frac{a_0}{r}$ occurs at $r = \frac{\sqrt{a_0}}{\sqrt{a_1} - \sqrt{a_0}}$. Therefore, dividing the numerator and denominator of \eqref{eqn:useful2} by
$$M=\frac{a_1}{1+\frac{\sqrt{a_0}}{\sqrt{a_1} - \sqrt{a_0}}} - \frac{a_0}{\frac{\sqrt{a_0}}{\sqrt{a_1} - \sqrt{a_0}}},$$
we see that
\begin{equation}\label{eqn:useful3}
    \beta\leq 1 - \frac{1}{\frac{c_1}{M}\log d_L + \frac{c_2}{M} n_L \log t + \frac{c_3}{M} n_L + \frac{c_4}{M}}.
\end{equation}
In Table \ref{tab:zero_free_regions_final}, we present the constants for two choices of $\varepsilon$. Observing the values for $\varepsilon = 0.01$, inequality \eqref{eqn:useful3} will yield the explicit zero-free region \eqref{eqn:main_result1} for $t\geq 1$, which completes the proof of Theorem \ref{thm:main_result1}.

\begin{table}[h!]
    \centering
    \begin{tabular}{l|ll}
        & $\varepsilon = 0.15$ & $\varepsilon = 0.01$\\
        \hline
        $M$              & 0.1021253857  & 0.1021253857\\
        $\frac{c_1}{M}$  & 12.24106100   & 12.24106100\\
        $\frac{c_2}{M}$   & 9.534650638   & 9.534650638\\
        $\frac{c_3}{M}$  & 0.444485082  & {0.050168175}\\
        $\frac{c_4}{M}$   & 5.123026304   & 2.269182727
    \end{tabular}
    \caption{Constants for the explicit zero-free region in Theorem \ref{thm:main_result1} given $\varepsilon = 0.15$ or $\varepsilon = 0.01$.}
    \label{tab:zero_free_regions_final}
\end{table}

\section{Proof of Theorem \ref{thm:main_result2}}

Theorem \ref{thm:main_result2} is an improvement of part of \cite[Theorem 1.2]{kadiri2012explicit}. We will recycle Kadiri's proof, except we use the polynomial $p_{16}(\varphi)$ in place of a polynomial from $P_4$. Suppose $\log d_L$ is asymptotically large and consider three regions,
$$\mathbb{I}_A = \left(0,\frac{d_1}{\log d_L}\right]\mbox{, }\mathbb{I}_B = \left(\frac{d_1}{\log d_L},\frac{d_2}{\log d_L}\right],\mathbb{I}_C = \left(\frac{d_2}{\log d_L}, 1\right),$$
where $d_1, d_2$ are constants to be chosen. Suppose further, that 
$$\sigma - 1 = \frac{r}{\log d_L}\mbox{ and } 1 - \beta = \frac{c}{\log d_L}.$$
In the regions $\mathbb{I}_B$ and $\mathbb{I}_C$, we impose further restrictions. Suppose $0 < c, r < 1$ such that
$$\frac{a_0}{a_1 - a_0}c < r \mbox{ and }d_2 > \frac{\sqrt{r(r+ c)}}{2}.$$

Combining analogous arguments to those results in \cite[§3.2, §3.3, §3.4]{kadiri2012explicit}, one can easily establish that
\begin{equation}\label{eqn:IA_kadiri}
    0\leq
    \frac{1}{r} - 2\frac{r+c}{(r+c)^2 + {d_1}^2} + \frac{1 - \kappa}{2}
\end{equation}
in the region $\mathbb{I}_A$,
\begin{align}
    0\leq &\mcE_B(d_1, d_2, r, c)\nonumber\\
    &:= \frac{a_0}{r} - \frac{a_1}{r+c}
    + \frac{a_1 r}{r^2 + {d_1}^2}
    - \frac{a_0(r+c)}{(r+c)^2 + {d_1}^2}\nonumber\\
    &-\frac{a_0(r+c)}{(r+c)^2 + {d_2}^2}
    -\frac{a_1(r+c)}{(r+c)^2 + 4 {d_2}^2}
    + \frac{1 - \kappa}{2}\sum_{k=0}^{16}a_k\label{eqn:IB_kadiri}\\
    &+\sum_{k=2}^{16} a_k\left(\frac{r}{r^2 + k^2{d_1}^2}
    - \frac{r + c}{(r+c)^2 + (k-1)^2{d_2}^2}
    - \frac{r + c}{(r+c)^2 + (k+1)^2{d_2}^2}\right)\nonumber
\end{align}
in the region $\mathbb{I}_B$ and
\begin{align}
    0&\leq \mcE_C(d_2, r, c)\nonumber\\
    &:= \frac{a_0}{r} - \frac{a_1}{r+c} + \frac{a_1 r}{r^2 + {d_2}^2} - \frac{a_0(r+c)}{(r+c)^2 + {d_2}^2} + \frac{1 - \kappa}{2}\sum_{k=0}^{16}a_k \label{eqn:IC_kadiri}
\end{align}
in the region $\mathbb{I}_C$. Suppose $d_1$ and $r$ are fixed. The admissible values of $c$ which one can input into \eqref{eqn:IA_kadiri} are those $c$ such that 
\begin{equation}\label{eqn:IA_solver}
    c\geq \frac{\sqrt{r^2 - {d_1}^2\left(1 + \frac{1 - \kappa}{2}r\right)^2} - \frac{1 - \kappa}{2} r^2}{1 + \frac{1 - \kappa}{2} r}.
\end{equation}
Denote the smallest value for $c$ in \eqref{eqn:IA_solver} by $c_A$. Next, let $c_B$ denote the root of $\mcE_B(d_1, d_2, r, c)$, where $r$ is chosen such that the root $c_B$ is as small as possible. Similarly, let $c_C$ denote the smallest root of $\mcE_C(d_2, r, c)$ for some optimally chosen $r$. It follows that $\zeta_L$ has at most one zero in the region $s=\sigma + it$ such that $t < 1$ and
$$\sigma \geq 1 - \frac{1}{R \log d_L}$$
such that $R = \max\left(1/{c_A}, 1/{c_B}, 1/{c_C}\right)$. Moreover, if an exceptional zero exists then it is real and simple by \cite[§3.5]{kadiri2012explicit}. To complete our proof of Theorem \ref{thm:main_result2}, it will suffice to show that $R = 12.43436$ is an admissible value.

First, suppose that we choose the same values that Kadiri chose; $d_1 = 1.021$ and $d_2 = 2.374$. One can establish that $1/{c_A} = 12.5494$ when $r = 2.1426$. Moreover, using our higher degree polynomial, we can compute the roots of $\mcE_B(1.021, 2.374, r, c)$ and $\mcE_C(2.374, r, c)$ over a selection of $r$. The results of these computations are presented below.
\begin{center}
\begin{tabular}{c|c|c}
     \textit{Root of} & $r$ & $1/c$\\
     \hline
     $\mcE_B(1.021, 2.374, r, c)$ & 0.2366 & 12.43922\\ 
     $\mcE_C(2.374, r, c)$ & 0.2477 & 12.42548
\end{tabular}
\end{center}
Therefore, these choices of $d_1$ and $d_2$ would yield Theorem \ref{thm:main_result2} with
$$R = \max\left(12.5494, 12.43922, 12.42548\right) = 12.5494.$$

Above, the limiting factor appears to be the value for $1/{c_A}$. We can reduce the value of $1/{c_A}$ by decreasing the value of $d_1$, however, we are also limited by the sizes of $1/{c_B}$ and $1/{c_C}$ which we can obtain. Therefore, we only need to choose $d_1$ such that $1/{c_A}$ is small enough. The cost of choosing $d_1$ too small is a larger interval $\mathbb{I}_B$, which might not be ideal.

Given $d_1$, to find a good enough choice for $d_2$, we have tested many values for $d_2$ and computed the optimal outcomes in each case. If one chooses $d_1 = 1.0015$, then we found (to 3 decimal places) that $d_2 = 2.318$ yielded the best results. For this $d_1$, one can determine that $1/{c_A} = 9.7946$ when $r = 2.1163$. The results of the remaining computations for $1/{c_B}$ and $1/{c_C}$ are presented below.
\begin{center}
\begin{tabular}{c|c|c}
     \textit{Root of} & $r$ & $\frac{1}{c}$\\
     \hline
     $\mcE_B(1.0015, 2.318, r, c)$ & 0.2363 & 12.43355\\ 
     $\mcE_C(2.318, r, c)$ & 0.2473 & 12.43436
\end{tabular}
\end{center}
Therefore --- as required --- these choices of $d_1$ and $d_2$ will yield Theorem \ref{thm:main_result2} with
$$R = \max\left(9.7946, 12.43355, 12.43436\right) = 12.43436.$$

\bibliographystyle{amsplain}
\bibliography{refs}

\end{document}